# Cold-Tip Temperature Control of Space-borne SatelliteStirlingCryocooler: Mathematical Modeling and Control Investigation

Shaival H. Nagarsheth*. Jiten H. Bhatt.**
Jayesh J. Barve***

* PG. Scholar, Instrumentation & Control, Nirma University, Ahmedabad, INDIA(E-mail: shn411@gmail.com)
(Currently: PhD Scholar, Electrical Engineering, SVNIT, Surat, INDIA.)
**External PhD Scholar,Nirma University, Ahmedabad, INDIA (E-mail: jiten@sac.isro.gov.in)
(Also: Scientist/Engineer, Space Applications Centre, Indian Space Research Organisation,ISRO, Ahmedabad,INDIA
***Former Professor, Nirma University, Ahmedabad, India (E-mail: barve.jayesh@ieee.org}
(Currently:Principal Engineer, Controls & Optimization,GE Global Research Centre, JFWTC, Bangalore, INDIA)

Abstract: IR (Infra-Red) detectors are widely used in Space-borne remote sensing satellites. In order to achieve a high signal to noise ratio, the IR detectors need to be operated at cryogenic temperatures. Traditionally, the cryogenic cooling of these detectors is achieved using passive cooling techniques. However recent trend is to employ Stirling-cycle based miniaturized active cryocoolers. An accurate and stringent control of active cryocooler cold-tip temperature is essential to accomplish high signal & image quality from the IR detectors. This paper presents work on investigations and comparison of performance of proposed 2-DOF (2-Degrees-of-Freedom) versus traditional 1-DOF feedback-control structures for the control of cryocooler cold-tip temperature used in IR (Infra-Red) detectors of Space Satellites. Towards this, first-principle based control oriented mathematical model simulated in Matlab/Simulink is proposed to support such investigation and controller tuning. Open-loop (system) and closed-loop (controls) simulation results are tuned & validated with the experimental data obtained from the Lab-scale test-setup of a commercial Stirlingcryocooler. The performance of 2-DOF feedback control structures are analyzed with the help of experiment results offering a better cold-tip temperature control performance.

*Keywords:* Cryogenic, Mathematical Modeling, 2DOF-PI, Temperature Control, Infrared detectors

## 1. INTRODUCTION

The Infrared (IR) detectors for space borne remote sensing applications need to operate at cryogenic temperatures since they get exposed to noise background. In the past, detector cooling was achieved using passive radiant coolers(Bard S. 1984). The limited capacity of passive coolers often limits the number of detectors and the detector power for space borne camera(Zhang et al. 2010). With the advancement of miniaturized Stirling cycle based active cryocoolers, the next generation IR remote sensing payloads employ large-area array detectors with a large power dissipation capability. The use of large area array detectors often demands stringent cold-tip temperature stability(Ross R. J. 2001) and regulation/control to a) obtain uniform temperature across IR detectors and b) maintain a uniform leakage current from all the IR detector-pixels. The stringent cold-tip temperature stability and control requirements pose a challenge to traditional cryocooler controllers, motivating to investigate and design more appropriate control structure and well-tuned well-performing controllers. There exists various advanced model based predictive multivariable and optimizing control algorithms that can be used potentially(Carlos E. G. et al. 1989)(Rafal N et al. 2015). But, on board satellite requires use of computational viable simpler PID control algorithms due to constraint of computational memory and computational time. In practice, mostly PI controllers are employed and considered sufficient for the cryocooler control electronics purpose.

Typically, cryocooler cold tip temperature control requirement is very stringent with permissible cold-tip temperature error being < ±0.1K. Also, the space-borne instruments employing IR detectors operates IR detectors and cryocoolers on the need basis and turns off the IR detectors and cryocoolers to conserve the Satellite power during the rest of the period. Each detector mode has got its own power requirements and in effect leads to load disturbances at the cold tip of cryocooler. For having a quick cool down and low power consumption, it is also planned to operate the cryocooler at higher cold tip temperature during non-imaging periods. Hence considering these different requirements, it is essential to have a cryocooler controller that provides a good set-point tracking as well as disturbance-rejection responses for the cryocooler cold-tip temperature. It is found that researchers have also carried outwork on dynamic simulation through mathematical modelling (Guo. D. et al. 2013)(Zhang et al. 2002) and on a scheme for stroke length control with robustness test on change of load assuming FOPDT transfer function model (Yang et al. 1999). However, detailed mathematical model with the consideration of heat and mass transfer inside the cryocooler, based on gas equation and thermodynamics is not cited. Also, detailed parametric study for load disturbance (continuous and pulsating) and change in input parameters has not been cited.

In this research work, the authors have attempted to carry out investigations by 1-DOF (degree of freedom) control structure and then by 2-DOF control-structure implemented practically for the cold-tip temperature control of acryocooler used for cooling IR detectors in space bornesatellite. Towards these investigations, first-principle energy-balance based mathematical model is proposed. The proposed mathematical model is simulated in Matlab-Simulink platform (S - function builder tool) to help controller tuning and validated with experimental results. The model is tuned using experimental data obtained from the lab-scale test setup of a commercial Stirling cryocooler (used with the Space Satellite IR detectors).

The paper is organized as follows. The section-1 presents motivation for this work with some relevant literature references. Section-2 discusses the background covering the proposed control-oriented transient mathematical model of the Stirlingcryocooler, 1-DOF & 2-DOF control structures. Section-3 provides a brief-overview of cryocooler control-electronics of the NI-LabVIEW based Lab-scale experimental test set-up. In Section-4, simulation and experimental test-results for (a) tuning & validation of the proposed mathematical model and (b) performance comparison & investigations of traditional 1-DOF versus proposed 2-DOF controller are presented and discussed. Section-5 presents the concluding remarks. Relevant references of published literature and acknowledgements are given at last.

## 2. BACKGROUND

*2.1 Mathematical Modelling*

Stirling cryocooler operation is considered as per the schematic shown in Fig.1 for modelling and further investigations. (Randall 1985)

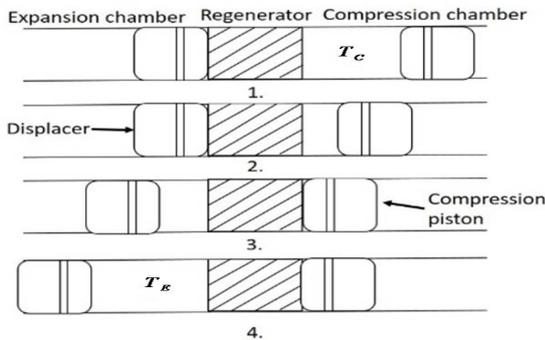

Fig. 1. Schematic of stirlingcryocooler Operation.

Fig.1 depicts a cylinder containing two opposed piston arrangement with a regenerator in between them. It has two volume spaces, a compression space at temperature $T_C(K)$ and an expansion space at temperature $T_E(K)$. The temperature gradient between the two chambers is $T_C - T_E$ assuming that there is no thermal conduction in the longitudinal direction. Both the pistons move without friction or leakage of the working fluid enclosed between them.
The compression piston is made to reciprocate with the help of a linear motor. Piston motion and the compression is considered as a mass spring damper system assembly (Yang et al. 1999) governed by the equation:

$$M\ddot{x} + D\dot{x} + Kx = \lambda IB - A\Delta P \qquad (1)$$

Where $\lambda IB - A\Delta P$ is the force to the compression piston, $M(kg)$ is mass of cylinder, $x$ is displacement of the piston, $D$ is damping co-efficient, $K$ is friction, $\lambda$ is length of coil, $I$ is input current to motor, $B$ is magnetic flux density, $\Delta P$ is the pressure difference between the two chamber.

Stirlingcryocooler is governed by the ideal Stirling cycle as shown in Fig. 2 represented on P-V diagram (Thomas F. 2005)

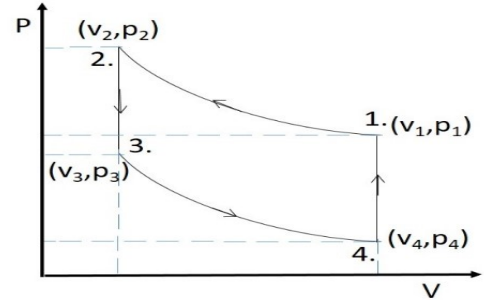

Fig. 2. Ideal Stirling cycle- PV diagram.

The ideal cycle consists of four processes as below:
Process (1-2) is compression, where piston moves from the outer dead point, at that time the expansion space piston is at the inner dead point and remains stationary. The working fluid (helium) in the compression space is at ambient temperature, volume and pressure are at state (1) on the P-V diagram in Fig. 2.
In the compression space helium gas is compressed so the pressure and temperature increases as the volume decreases. The heat generated is rejected to the sink at rate ($\dot{Q}_{rej}$). The governing equations to this process are:

According to ideal gas law:
$$PV = nRT$$
$$\therefore PdV + VdP = nRdT$$
Rate of work done is given by
$$P\,dV/dt + V\,dP/dt = nR\,dT/dt = \dot{W}_{1-2} \qquad (2)$$
Where $\dot{W}_{1-2}$ is the rate of change of work done on the system, $R$ is gas constant, $n$ number of mole of gas and the changes in the volume and pressure will lead to state (2) on P-V diagram.

According to first law of thermodynamics, change in the internal energy (heating) equal to work done and in addition heat rejection(Incropera F. P 1990).

$$\therefore M_{he}C_{phe}\left(\frac{dT_C}{dt}\right) = \dot{W}_{1-2} - \dot{Q}_{rej} \qquad (3)$$

Where $M_{he}(kg)$ is mass of helium gas in the compression chamber which is compressed, $C_{phe}(J/kg\,K)$ is specific heat of helium gas assumed to be constant with change in temperature, $T_C$ is the temperature of the working fluid (helium) in the compression chamber whose change leads to change in the internal energy and $\dot{Q}_{rej}(J/s)$ is the rate of heat released to the sink.

Process (2-3) is a constant volume heat transfer process where both pistons move simultaneously. Working fluid passes through the regenerator and is cooled from $T_C$ to $T_E$ by giving heat $Q_c(J)$ to the regenerator matrix (Sercan 2011). Passing through matrix at constant volume causes decreases in pressure.

$$Q_c = \varepsilon\, M_{min} C_{min}(T_C - T_{R1}), \qquad (4)$$
$$Q_C = M_r C_{pr}(T_{R1} - T_{R2}), \qquad (5)$$
$$Q_C = M_{he} C_{phe}(T_C - T_E). \qquad (6)$$

Where, $M_r(kg)$ is the mass of regenerator material, $\varepsilon$ is the effectiveness of the regenerator which can be obtained by (e-NTU) number of transfer unit method (C. J. Paul et al. 2015) (Urieli I et al. 1984). $T_{R1}(K)$ and $T_{R2}(K)$ are the temperatures of regenerator before and after the Process(2-3) respectively. $C_{pr}(J/kg\,K)$ is the specific heat of regenerator material.

During process (3-4) expansion takes place as the expansion chamber piston moves away from the regenerator providing useful refrigeration of the cycle (cooling). So here at the cold tip in the expansion chamber where the effective cooling is achieved, heat at rate $\dot{Q}_{ab}(J/s)$ is added at this juncture from an external heat source (load) to the working fluid which can be described by the following equation:

$$M_{he} C_{phe}\left(\frac{dT_E}{dt}\right) = \dot{W}_{3-4} + \dot{Q}_{ab} + \dot{Q}_{bl} \qquad (7)$$

Where $\dot{Q}_{bl} = hA_s(T_{ab} - T_E)$ is the rate at which heat is added from the base load (surrounding) which is at ambient temperature $T_{ab}(K)$.

Process (4-1) is again a constant volume heat exchange but over here the regenerator gives away all the heat $Q'_c(J)$ consumed in process (2-3) again to working fluid (helium) such that volume, pressure and temperature of the system returns to the state (1) shown in Fig. 2.

$$Q'_c = \varepsilon\, M_{min} C_{min}(T_{R2} - T_E), \qquad (8)$$
$$Q'_C = M_r C_{pr}(T_{R2} - T_{R1}), \qquad (9)$$
$$Q'_C = M_{he} C_{phe}(T_E - T_C). \qquad (10)$$

*2.2 Traditional 1-DOF Controller*

Fig. 3 shows a conventional one degree of freedom controller (1DOF). In general, such 1-DOF PI controller is implemented using a classical PI controller for the cryocooler temperature control.

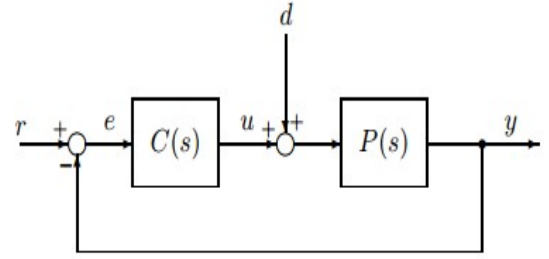

Fig. 3. Conventional 1-DOF Feedback Controller

However, such 1-DOF structure poses a challenge in accomplishing high accuracy (particularly for tight temperature control), high performance for set-point tracking and disturbance-rejection simultaneously (M. Araki et al. 2003). As first attempt in the present work the above shown conventional 1-DOF feedback controller is implemented for the cold tip temperature control to confirm the predictions of M. Araki.

*2.3 Proposed 2-DOF controller*

Two degree of freedom (2DOF) control system has advantages over a 1DOF system (I. M. Horowitz 1963). The extra parameter that the 2DOF control algorithm provides is used to improve their servo-control behavior while considering the regulatory control performance and the closed-loop control system robustness. Considering these advantages of 2DOF control, a digital PI controller with 2DOF has been implemented for the cold tip temperature control of cryocooler.
Many different equivalent forms of 2DOF exists (M. Araki et al. 2003)(V. M. Alfaro et al. 2016). The most suitable form for implementation for cryocooler control is the Set-point filter type expression of the 2DOF PID control system which is shown in Fig. 4.

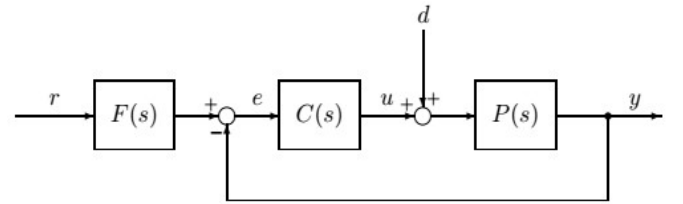

Fig.4. Block Diagram of Set-point filter 2-DOF Controller

The transfer function (Goodwin 2001) of the Set Point Filter and Controller for Filter type of 2DOF PID controller is given by:

$$F(s) = \frac{1}{1 + \tau_f s} \qquad (11)$$

$$C(s) = K_p\left(1 + \frac{1}{\tau_i s}\right) \quad (12)$$

This strategy introduces changes in set-point gradually rather than abruptly. The filter time constant $(\tau_f)$ determines how fast the filtered set-point will attain the new value. This can eliminate or can effectively reduce overshoot for set-point (Seborg et al. 2004). Because of the smooth set point response, the controller can now accommodate higher Proportional gain $(K_p)$ to improve the disturbance-rejection performance without compromising on the set-point tracking performance.

### 3. CRYOCOOLER CONTROL TEST SET-UP

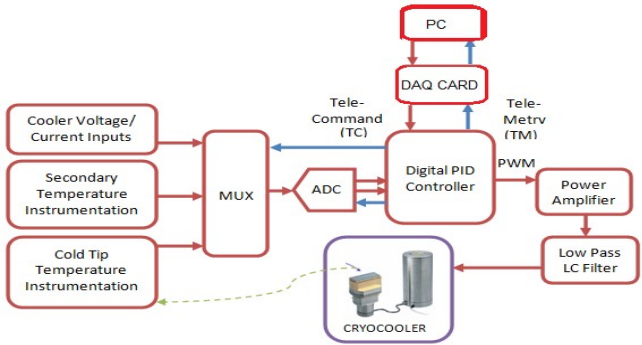

Fig. 5. Block Diagram of Cryocooler Control Electronics

Fig.5 shows a simplified block diagram of Cryocooler Controls Electronics. A digital PI controller controls the PWM output signal which is given to a Class D power amplifier which drives the coils of cryocooler. The controller also implements auxiliary functions such as operating mode selection using Tele-Command (TC) and health parameter monitoring using Telemetry (TM). The TC and TM interface is simulated by a PCI based DAQ card. NI LabVIEW based software is used to issue Tele-commands and receive the Telemetry. The cold tip temperature is sensed by a diode based temperature sensor. The linear drive Stirling cryocoolers need an AC signal drive, for which Direct Digital Synthesis (DDS) technique is used in Digital Controller. In this technique, the sine wave coefficients stored in a ROM are used to generate a digital sine signal which is then transmitted to the power amplifier using PWM techniqueFig. 6 shows the photograph of the complete test set up.

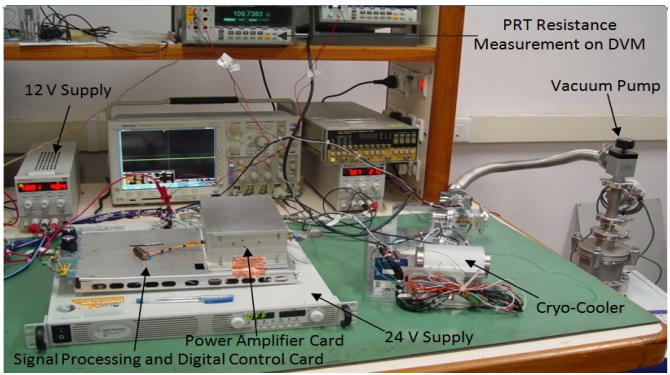

Fig. 6. Cryocooler Test set up &it's Control Electronics

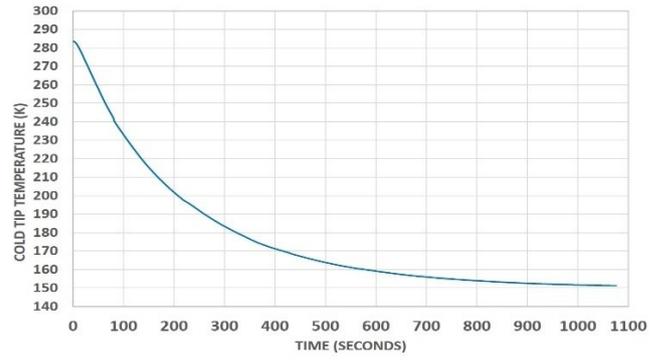

Fig. 7. Start-up cool-down: Cold-tip temperature response for 7.27 Watts cryocooler input power

Fig. 7 shows the results of the cool down of the Stirlingcryocooler obtained practically with initial power of 7.27W. Parametric study is carried out on the attainment of steady state temperature after cool down. Here for cool down the cold tip temperature $(T_E)$ has been taken into consideration and the cool down time is found out to be 20 minutes.

### 3. RESULTS AND DISCUSSION

Implementation of the mathematical model described from (1) to (10) is carried out using S-function in Matlab/Simulink platform to perform analysis. For validation of model the tested machine is a linear Stirling cryocooler and the cycle is assumed to start at room temperature (300 K) position (1) in PV diagram. The main parameters are as follows: cylinder mass is 1kg, damping co-efficient is 0.116 N-Sec/m, spring constant is 40000 N/m. Air-gapflux density is 0.25T, operating pressure 20 bar, length of compressor 0.06m, length of expander is 0.082m, diameter of expander is 0.013m,effectiveness of regenerator is 0.98, free convective heat transfer coefficient $h$ for base load heat transfer is $50.481 \, W/m^2 K$.

The parametric study is carried out to observe the cold tip temperature variation against step-change in
(a) The compressor motor input current $-I$
(b) The cooling load rate $-\dot{Q}_{ab}$.
Simulator also records the time responses of various process variables like displacement, pressure, volume compressed, work done, and hot-tip temperature after every operational cycle (50Hz).

*3.1 Comparison of simulated and experimental temperature response for a step-change in the input motorcurrent*

Fig. 8 to Fig. 10 shows the comparison of simulated and experimental cold-tip temperature time-responses for three different step-signals applied to a compressor motor current $(I)$ around the initial steady-state operating condition of 151K cold-tip temperature.
Table-1 shows the steady-state values of cold-tip temperature for the same three step-changes in motor current. The time-response and steady-statevalues of model-simulation results

match quite well with the experimental results as observed in Fig 8-10 and Table-1.

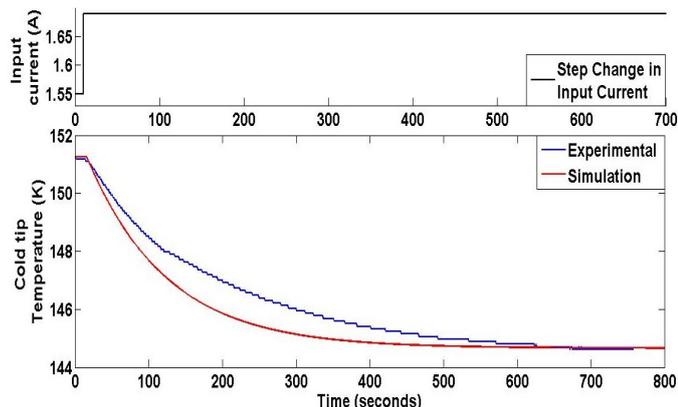

Fig. 8. Open loop response for step change in input change from 1.55A to 1.69A

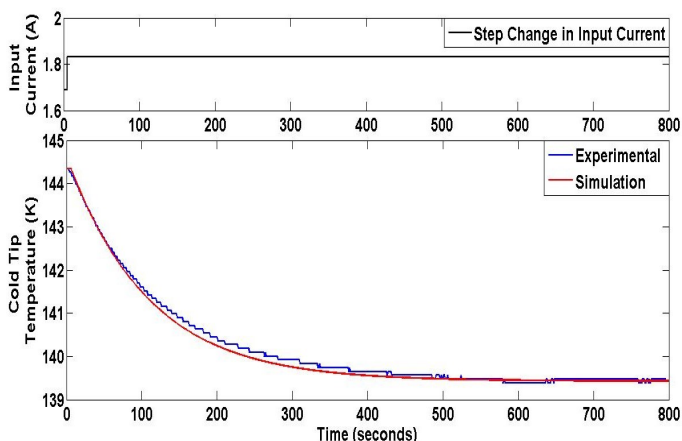

Fig. 9. Open loop response for step change in input change from 1.69A to 1.83A

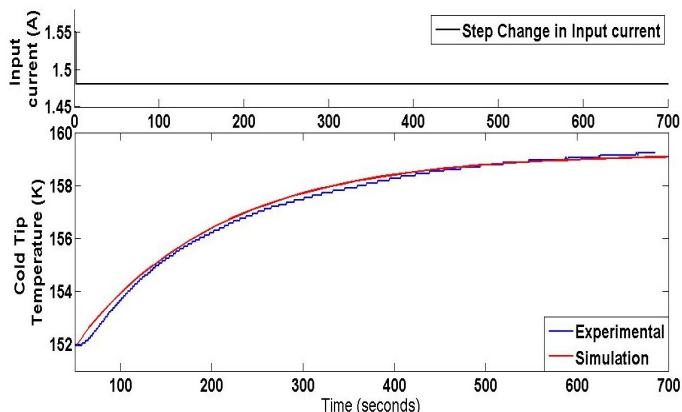

Fig. 10. Open loop response for step change in input change from 1.55A to 1.48A

**Table 1. Model validation: Response to motor current**

| Current I (A) | Cold Tip Temperature (K) | |
|---|---|---|
| | Experimental | Simulated |
| 1.55 | 151.26 | 151.3 |
| 1.69 | 144.66 | 144.7 |
| 1.83 | 139.47 | 139.5 |
| 1.48 | 159.23 | 159.1 |

*3.2 Comparison of simulated and experimental temperature response for a step-change in the cooling-load (disturbance).*

Similarly, Fig. 11 to Fig. 13 shows the comparison of simulated and experimental cold-tip temperature time-responses for 3 different pulse-signals applied to cryocoolerload around the initial steady-state operating condition of 151K cold-tip temperature. Table-2 shows thevalues of main parameters of the cold-tip temperature time-response for the same three pulse-signals applied to cryocooler cooling-load. Here, the time-response and key time-response parameters of the model-simulation results match quite well with the experimental results as observed in Fig. 11 to Fig. 13 and Table-2.

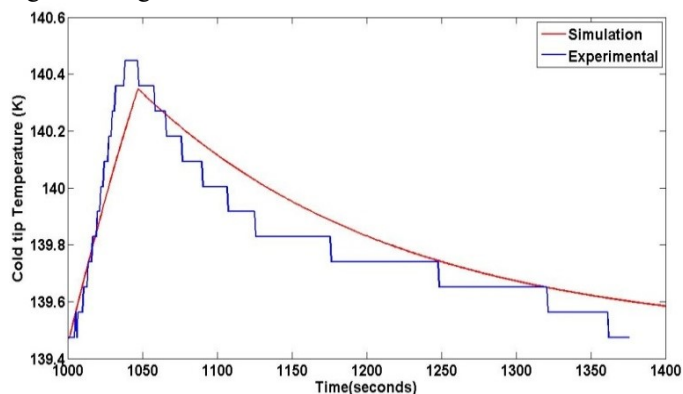

Fig. 11. Cold-tip temperature response for 118.8mW, 45 sec pulse-load applied to cooler load

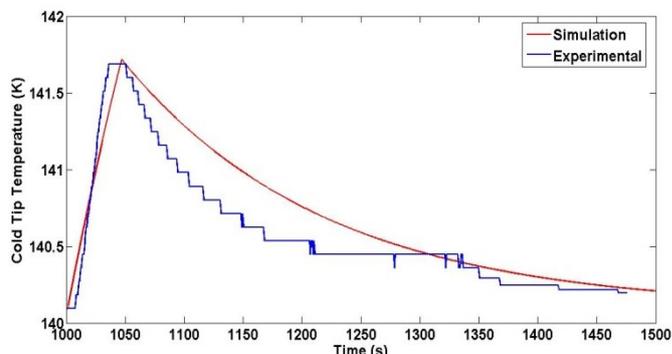

Fig. 12. Cold-tip temperature response for 220.5mW, 45 sec pulse-load applied to cooler load

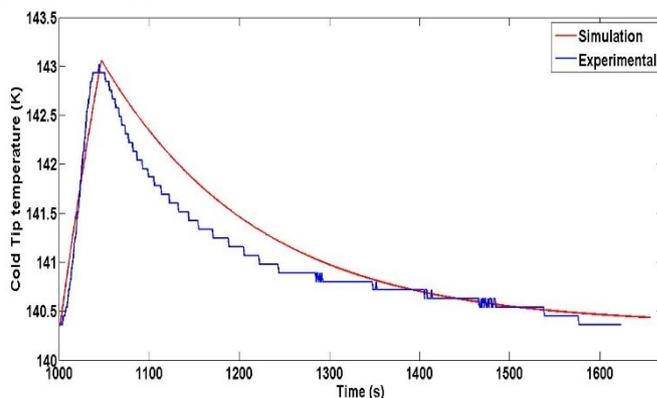

Fig. 13. Cold-tip temperature response for 330mW, 45 sec pulse-load applied to cooler load

Table-2. Model validation: Load-change Response

| Thermal Load-Pulse(45 s) | Peak Temperature (K) | | Difference in peak value |
|---|---|---|---|
| | Simulation values | Experimental values | |
| 118.8mW | 140.30 | 140.4 | 0.10 |
| 220.5mW | 141.61 | 141.69 | 0.08 |
| 330mW | 142.90 | 143.1 | 0.20 |

*3.3 Investigations of 1-DOF and 2-DOF PI controllers*

Next the closed loop control results for the set-point and disturbance-rejection responses of the same cryocooler cold-tip temperature are presented to investigate performancecomparison of traditional 1-DOF and proposed 2-DOF control schemes for the considered cold-tip temperature control system for the Satellite IR detector cryocooler system.

**1-DOF (One Degree of Freedom) Controller**
The initial tuning of PI Controller for cryocooler is based on technique described by(J Bhatt et al.2016), where in minimum number of hardware iterations are needed. Fig. 14 and Fig. 15 shows the experimental results for the set-point tracking and disturbance-rejection responses with two different sets of PI controller tuning parameters. If the controller parameters proportional-gain ($K_p$) andintegral-gain $K_i = K_p/\tau_i$ are tuned as (7.5 and 0.3), there is no peak-overshoot in theset-point tracking response. But, this controller now shows 1K peak-overshoot for even smaller 0.5W of load disturbance. It is further observed that if a proportional-gain is doubled (without change in the integral-gain), the disturbance rejection improves as the peak-value of the cold-tip temperature reduces from 1K to 0.5K for the same load pulse. However, for this case, the set-point tracking response deteriorates showing a significant overshoot of ~3.4K. These results validate the statement that 1-DOF control structure cannot simultaneously meet set-point tracking and disturbance-rejection responses specifications due to lack of capability to tune or shape set-point tracking and disturbance-rejection performances simultaneously as given in (M. Araki et al.2013). Hence, proposed 2DOF PI controller is investigated in this study.

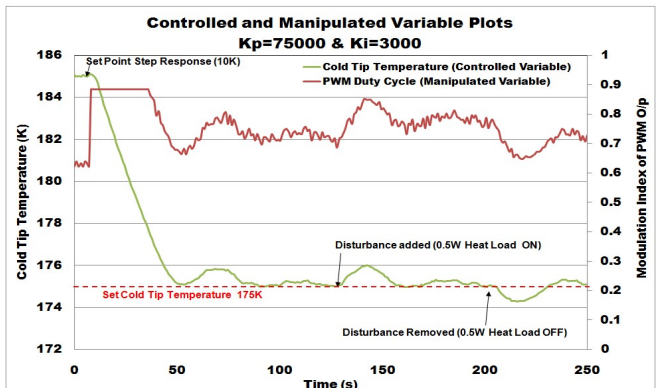

Fig.14. Set Point and Disturbance Response with 1DOF Digital PI Controller with Kp=7.5 & Ki=0.3

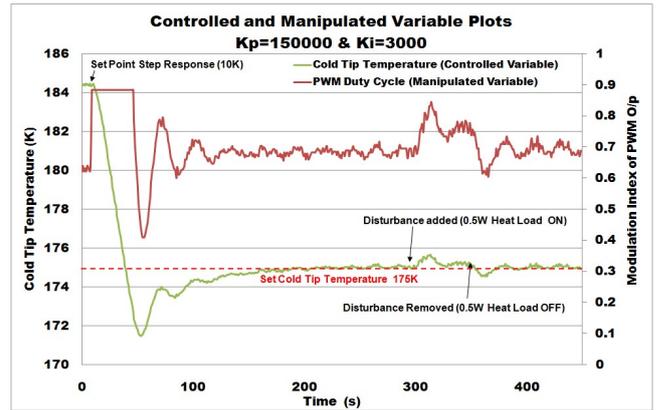

Fig. 15. Set Point and Disturbance Response with 1DOF Digital PI Controller with Kp=15 & Ki=0.3

**2-DOF (Two Degrees of Freedom) Controller**
To implement and evaluate 2-DOF control structure, set-point response filter block $F(s)$ shown in Fig. 4 is implemented as a first order recursive digital filter for set-point input signal which is implanted using the discrete equation given by:

$$y[n] = ay[n-1] + (1-a)x[n] \qquad (13)$$

Where, $y[n]$ is the output, $x[n]$ is the input and '$a$' is the filter-tuning parameter that is tuned based on the time-constant of the desired first-order step-response for the set-point tracking. In our investigations, value of '$a$' is tuned to 0.98 to obtain the effective time-constant of ~180 sec.

Fig. 16 shows the results obtained with two degree of freedom controller. As seen from the results it is now possible to increase proportional gain without getting any peak overshoot. By marginally increasing the proportional gain to 10 (while keeping the same integral gain of 0.3) the peak for load disturbance heat load of 0.5W reduces to only 0.5K. Thus it can be seen that it is now possible to optimize both set-point tracking response and disturbance rejection simultaneously.

Table-3 gives the comparison of various parameters for three different cases i.e. two cases representing two 1-DOF controllers (two sets of controller-parameters) and the remaining 1-case representing results of 2-DOF controller.

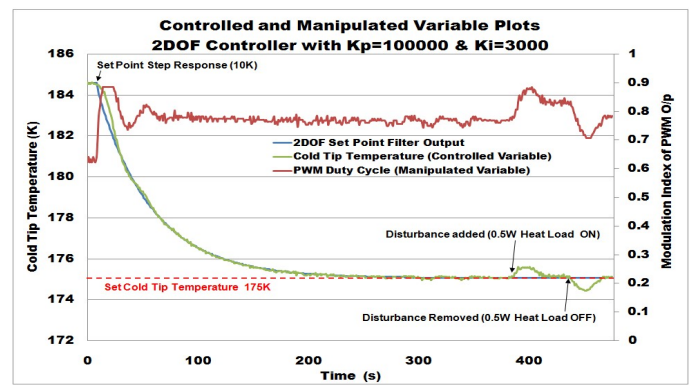

Fig. 16. Set Point and Disturbance Response with 2DOF Digital PI Controller with Kp=10 & Ki=0.3

**Table 3: Comparison: 1-DOF and 2-DOF controllers**

| Condition | Parameters | 1DOF<br>Kp= 7.5<br>Ki=0.3 | 1DOF<br>Kp=15<br>Ki=0.3 | 2DOF<br>Kp=10<br>Ki=0.3 |
|---|---|---|---|---|
| Set Point Response | Peak Time (s) | 43 | 48 | - |
|  | Peak Overshoot | - | 3.4 | - |
|  | Settling Time (s) | 72 | 207 | 240 |
| Disturbance Response | Peak Time (s) | 19 | 17 | 18 |
|  | Peak Overshoot | 1 | 0.5 | 0.5 |
|  | Settling Time (s) | 33 | 27 | 30 |

## 4. CONCLUSION

A mathematical model of the specific Stirling cryocooler is developed based on the first-principles i.e. thermodynamics and gas-law equations. The model is simulated in Matlab-Simulink platform. Several simulations are carried out to tune and validate the operational responses of the cryocooler. Besides, the cold-tip temperature response is simulated for the step-change applied to compressor-motor-current input and a pulse-change applied to thermal-load (disturbance). The simulation results are compared with the experimental results obtained by applying the similar step and pulse-change signals applied on the Lab-scale IR-detector Stirling Cryocooler test-setup. The model simulation results are found to match very well within ~0.2K for various thermal-load-change cases, and within ~0.1K for motor input current step-change. Traditional 1-DOF and proposed 2-DOF PI controllers are also investigated for the cold-tip temperature control. It is found that the proposed 2-DOF PI controller showsa superior performance compared to 1-DOF PI controllerallowing to simultaneously obtain good set-point tracking and disturbance-rejection responses. Next, we propose to extend this workby using this model & lab-setup platform to investigate other advanced controls methods, and also by developing multi-variable model for this cryo-cooler.

## ACKNOWLEDGEMENT

The authors would like to thank Mr. Manish Mehta and Mr. Nitin Upadhyay for their guidance and support in the experimental set-up and testing. Authors also acknowledge Mr. Sanjeev Mehta, Deputy Director, Sensor Electronics Development Area, and Mr. Arup Roy Chowdhury, Director, Space Applications Centre (SAC-ISRO) for providing an opportunity to work on this project.